\documentclass[final]{siamart171218}

\usepackage[T1]{fontenc}
\usepackage{mathtools}
\usepackage{amsfonts}
\usepackage{graphicx}
\usepackage{epstopdf}
\usepackage{algorithmic}
\usepackage{bm}
\usepackage{etoolbox}
\usepackage[caption=false]{subfig}
\usepackage{booktabs}
\usepackage{multirow}

\ifpdf
  \DeclareGraphicsExtensions{.eps,.pdf,.png,.jpg}
\else
  \DeclareGraphicsExtensions{.eps}
\fi

\newcommand{\MATLAB}{\textsc{MATLAB}\textsuperscript{\tiny\textregistered}}

\newsiamremark{remark}{Remark}
\newsiamremark{hypothesis}{Hypothesis}
\crefname{hypothesis}{Hypothesis}{Hypotheses}
\newsiamthm{claim}{Claim}

\headers{A Data-Driven Modeling Framework for Anomalous Rheology}{J.L. Suzuki, T.G. Tuttle, S. Roccabianca, M. Zayernouri}

\title{A Data-Driven Memory-Dependent Modeling Framework for Anomalous Rheology: Application to Urinary Bladder Tissue}

\author{Jorge L. Suzuki\thanks{Department of Mechanical Engineering and Computational Mathematics, Science and Engineering, Michigan State University
    (\email{suzukijo@egr.msu.edu}).}
\and Tyler G. Tuttle\thanks{Department of Mechanical Engineering, Michigan State University
    (\email{tuttlety@msu.edu}).} 
\and Sara Roccabianca\thanks{Department of Mechanical Engineering, Michigan State University
    (\email{roccabis@msu.edu}).}
  \and
  Mohsen Zayernouri\thanks{Department of Mechanical Engineering and Department of Statistics and Probability, Michigan State University (\email{zayern@egr.msu.edu}).}
  }

\usepackage{amsopn}

\makeatletter
\newcommand*{\addFileDependency}[1]{
  \typeout{(#1)}
  \@addtofilelist{#1}
  \IfFileExists{#1}{}{\typeout{No file #1.}}
}
\makeatother

\newcommand*{\myexternaldocument}[1]{%
    \externaldocument{#1}%
    \addFileDependency{#1.tex}%
    \addFileDependency{#1.aux}%
}

\ifpdf
\hypersetup{
  pdftitle={A Data-Driven Memory-Dependent Modeling Framework for Anomalous Rheology: Application to Urinary Bladder Tissue},
  pdfauthor={J.L. Suzuki, T.G. Tuttle, S. Roccabianca, M. Zayernouri}
}
\fi

\myexternaldocument{ex_supplement}

\begin{document}

\maketitle

\begin{abstract}
We introduce a data-driven fractional modeling framework for linear and nonlinear viscoelasticity aimed at complex materials, and particularly bio-tissues. From consecutive uniaxial multi-step relaxation experiments (up to $200\%$ strains) of five distinct anatomical locations of {porcine urinary bladder, we identify an anomalous relaxation character, with two power-law-like behaviors for short and long times, and strain nonlinearity for applied step strains greater than $25\%$}. Such behavior is affected by heterogeneous regional material compositions and nonlinearities taking place due to large strains. {The first component of our modeling framework is an \textit{existence study}, to determine admissible candidate fractional linear viscoelastic models, that qualitatively describe the (linear) time-dependent relaxation behavior. After the appropriate linear viscoelastic model is selected, the second stage include large-strain effects to the framework. For this purpose, among existing approaches in the literature, we utilize a fractional quasi-linear viscoelastic approach, given by a multiplicative kernel decomposition of the relaxation function from the selected linear model and a nonlinear elastic response for the bio-tissue of interest.} {From single-relaxation data of the urinary bladder, a fractional linear Maxwell model captures both short/long-term behaviors with two fractional orders, therefore being the most suitable candidate model for small strains at the first stage. For the second stage, multi-step relaxation data under large strains were employed to calibrate a four-parameter fractional quasi-linear viscoelastic model, that combines a Scott-Blair relaxation function and an exponential instantaneous stress response, the latter form being sufficient to describe the elastin/collagen phases of bladder rheology.} Our obtained results demonstrate that the employed fractional quasi-linear model, with a single fractional order in the range $\alpha = 0.25 - 0.30$ for distinct bladder samples, is suitable for the porcine urinary bladder, producing root mean squared errors below $2\%$ without need for recalibration over subsequent applied step strains. Our analyzes demonstrate that fractional models arise as attractive tools to capture the bladder tissue behavior under small-to-large strains and multiple time scales, therefore being potential alternatives to describe multiple stages of bladder functionality.
\end{abstract}



\section{Introduction}

Bio-tissues are complex and multi-functional materials, optimized for their specific host organisms, and constrained by limited set of building blocks and available resources \cite{Imbeni2005}. While the mechanical behavior of a number of standard engineering materials is quite well-understood, there is still a significant effort towards bio-materials, where microstructure heterogeneities, randomness and small scale physical mechanisms lead to non-standard and at times counter-intuitive responses. Power-law viscoelastic rheology is a complex response observed in many bio-tissues such as arteries \cite{Craiem2008}, cartilage \cite{Magin2010a}, lungs \cite{Suki1994}, smooth muscle \cite{Djordjevic2003}, liver and kidneys \cite{Nicolle2010}, among other classes of materials. These power-law materials, also termed \textit{anomalous}, exhibit {one or more power-law scalings} for creep/relaxation in the form $J(t) \propto t^{\beta}$ and $G(t) \propto t^{-\beta}$ {across multiple time-scales. Similar anomalous behaviors are also present for} dynamic storage/dissipation in the frequency domain \cite{McKinley2013,suzuki2021anomalous}. The origin of this power-law behavior at the continuum level is {linked} to (non-Fickian) sub-diffusive processes \cite{Metzler2000} in the corresponding fractal-like micro-structures \cite{Amblard1996}.

The aforementioned anomalous non-exponential behavior usually requires a {significant} number of material parameters when employing standard viscoelastic models. {These} consist of mechanical arrangements of linear springs and Newtonian dashpots, which induces a finite number of relaxation modes, which may lack predictability when performing outside the experimental {time scales \cite{Jaishankar2013,Magin2010a}. In this regard, fractional models become attractive alternatives, since their integro-differential operators naturally utilize power-law convolution kernels, coding self-similar microstructural features in a reduced-order mathematical language with smaller parameter spaces. Therefore, they have been employed as compact and predictive models for a number of anomalous systems, such as biological materials \cite{Magin2010a,Craiem2008,Suki1994,Djordjevic2003,Nicolle2010}, fluid turbulence \cite{samiee2020fractional,samiee2021tempered,akhavan2021data}, and instabilities \cite{akhavan2020anomalous}. We particularly note that such predictability has been shown to extend across different experiments (relaxation/creep) in certain cases \cite{Jaishankar2013}. Additionally, calibrating} experimental data with a set of existing rheological models leads to a material model selection problem, which is inherently ill-conditioned, since multiple models can pragmatically yield similar errors when confronted to experiment. In this work, we attempt to reduce this implicit ill-posedness by introducing fractional-order models as attractive alternatives to their integer-order counterparts, and employed to urinary bladder (UB) tissue modeling. Our fractional modeling framework aims to obtain compact mathematical models with a reduced number of material model parameters, while introducing a minimal, but sufficient number of fractional rheological elements that capture the qualitative response of multiple power laws and minimizes the errors, also rigorously taking into account the corresponding power-law memory effects.

The lower urinary tract, and especially the {UB}, is a highly dynamic organ system. To ensure its proper function, the bladder needs to be able to significantly increase in size while maintaining a low internal pressure, and this ability is dictated by the mechanics of the bladder wall. Specifically, during filling, the bladder tissue must leverage its viscoelastic characteristics to accommodate for large deformations without resulting in significant increase {of} luminal pressure. When this behavior is {compromised} due to disease, the resulting increase in pressure might generate a high-pressure urine reflux from the bladder to the kidneys, resulting in renal failure \cite{ansari2010risk,lopez2002bladder}. To increase the complexity of the organ mechanics, the characteristics of the bladder differ between different anatomical locations (i.e., dorsal, ventral, lateral, lower-body, trigone) \cite{korossis2009regional,MORALESORCAJO2018263} and orientations (i.e., longitudinal/apex-to-base and circumferential/transverse) \cite{korossis2009regional,MORALESORCAJO2018263,chen2013murine,gilbert2008collagen,cheng2018layer,van1978passive,NATALI20153088,wognum2009mechanical,nagle2017quantification}. To describe the mechanical behavior of bladder tissue, both hyperelastic \cite{cheng2018layer,wognum2009mechanical,nagle2017quantification,gilbert2008collagen,korkmaz2007simple,damaser1995effect,van1979contractility,van2002mechanical,regnier1983elastic,damaser1996two,watanabe1981finite,habteyes2017modeling} and viscoelastic \cite{van1978passive,NATALI20153088,coolsaet1973step,van1985model,coolsaet1975visco,van1981first,glerum1990mechanical,van1981dependence,alexander1973viscoplasticity,susset1981viscoelastic,kondo1973physical,alexander1971mechanical,venegas1991viscoelasticI,venegas1991viscoelasticII,nagatomi2004changes,nagatomi2008contribution} models have been used in the literature. {However,} due to the differences in mechanical testing protocols as well as modeling, most of the results cannot be compared with one another and often results in contradicting conclusions. While several pathologies of the lower urinary tract are associated with dramatic changes of the mechanical behavior of the bladder wall \cite{tuttle2021investigation}, still much is unknown about the mechanisms that affect this organ, not just in diseased states but in healthy as well. In this study, we focus on the healthy behavior of the porcine urinary bladder, which a present work suggested is a good model for the human urinary bladder.

Although fractional linear viscoelasticity has been succesfully employed in a number of biomechanical applications, additional modeling considerations are necessary when dealing with other material nonlinearities, such as large strains. This would imply that the material relaxation behavior depends on both time and applied strains, which requires additional modeling considerations. {In the UB case, Korossis et al. \cite{korossis2009regional} reports the small-strain regime to be drive by elastin, while a stiffer response for large strains was driven by collagen. Jokandan et al. \cite{JOKANDAN201892} characterized the quasi-static stress-strain response of the UB as exponential-like. To address such behavior,} a practical and {well-known class of models is the quasi-linear-viscoelastic (QLV) theory by Fung \cite{fung2013biomechanics}}, which considers a multiplicative coupling between a linear viscoelastic relaxation and a nonlinear elasticity term. {While the original formulation utilizes a finite spectrum of relaxation times, fractional extensions of the QLV theory have been developed and employed for the response of bio-tissues \cite{doehring2005fractional,Craiem2008}. Doehring et al. \cite{doehring2005fractional} employed a Mittag-Leffler-type reduced relaxation function that captured the short and long term behaviors of aortic valve cusp. Craiem et al. \cite{Craiem2008} developed a fractional Kelvin-Voigt-type reduced relaxation function with an exponential instantaneous response, which was successfully applied to the nonlinear viscoelasticity of arterial walls. Regarding the nonlinear viscoelastic modeling of the UB, Natali et al. \cite{NATALI20153088}, developed an anisotropic, visco-hyperelastic model, which was validated through a uniaxial experiment under cyclic loads. Nagatomi et al. \cite{nagatomi2004changes} studied the nonlinear behavior of rat bladders, by calibrating a two-dimensional QLV model to bi-axial relaxation data from bladders of subjects that were healthy and with spinal cord injury. Their findings reported a need for new models to account for both normal and pathological states, due to tissue remodeling.}

To the authors' best understanding, although existing studies have addressed the {nonlinear viscoelastic of the UB for different subjects and loading conditions,} there are no studies in the literature {leveraging the use of fractional viscoelastic models to model potential, emerging power-law behaviors}. In this work we develop a data-driven fractional modeling framework for linear and quasi-linear viscoelasticity to account for {both} anomalous power-law relaxation and large strains {of bio-tissues}. We validate the developed framework for the first time in the uniaxial relaxation of porcine urinary bladder tissue for a wide range of applied strains. The characteristics of our experimental procedure are:
\begin{itemize}
\item We obtain the porcine {UB} uniaxial relaxation data from small-to-large strains of five distinct anatomical locations.

\item Our relaxation experiments are performed under increasingly larger strains, without intermediate unloading steps or tissue preconditioning.

\item The mechanical response of the UB indicates nonlinear viscoelastic behavior with power-law relaxation, characterizing an anomalous{, non-exponential} behavior.

\end{itemize}
Given the {anomalous nonlinear response of the UB tissue, we develop our two-stage} anomalous modeling framework as follows:
\begin{itemize}
\item {In the first stage, we develop an \textit{existence} study that considers a set of linear fractional building block models (Scott-Blair, fractional Kelvin-Voigt, fractional Maxwell), which are selected according to the multi-power-law nature of the relaxation data and calibration errors at the linear viscoelastic regime.}

\item {In the second stage, we account for the large strain behavior of the corresponding tissue, by employing a fractional quasi-linear viscoelastic (FQLV). The goal is to extend the quality of power-law relaxation (due to the material's fractal microstructure) to the large strain regime of the tissue of interest.}
\end{itemize}

{We employed the aforementioned two-stage formulation to the UB experimental data, and obtained the following main findings;}

\begin{itemize}
    \item {All candidate} fractional linear viscoelastic models provide sufficiently accurate fits for single-relaxation steps under smaller strain levels, where the {two fractional order Maxwell model is the most suited for the UB data.}

    \item {The employed four-parameter FQLV model with a reduced Scott-Blair relaxation function and exponential instantaneous stress response was successful over five consecutive relaxation steps, with root mean squared errors below $2\%$, and without the need of model recalibration between applied step strains.}

    \item {The lower-range of obtained fractional orders is around $\alpha = 0.17-0.30$, which is compatible with the observed long-time slopes of the UB relaxation data. Small $\alpha$ values have been suggested to indicate strong fractality in bio-tissue microstructures such as collagen fibers \cite{doehring2005fractional}.}
\end{itemize}

The rest of the paper is organized as follows. In Section \ref{sec:setupmethod} we present the problem setup and methodology, with the uniaxial UB stress relaxation experiments and our proposed fractional modeling framework for biotissues, comprised of linear and quasi-linear fractional models. In Section \ref{sec:results}, we present our obtained linear viscoelasticity results for the UB relaxation under the first strain step, and the fractional quasi-linear viscoelastic model for all consecutive strain steps, followed by the discussions and potential improvements of the work. Finally, in Section \ref{Sec:Conclusions}, we provide the concluding remarks of this work and future directions.


\section{\label{sec:setupmethod}Problem Setup and Methodology}

\subsection{\label{sec:UBtests}Urinary bladder experimental relaxation tests}

\begin{figure}[ht!]
  \includegraphics[width=\columnwidth]{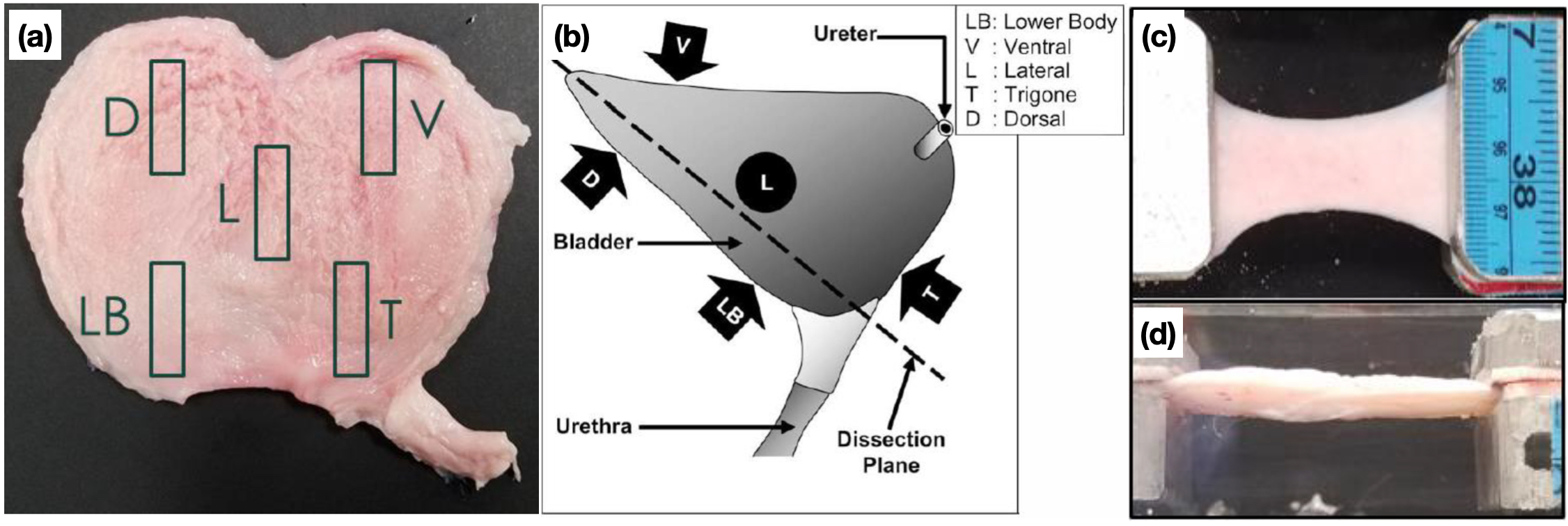}
	\caption{\label{fig:experimental}\textit{(a)} Dissected porcine UB showing the distinct anatomical locations from which samples were punched. \textit{(b)} A diagram of the UB from the lateral view. Source: \cite{korossis2009regional} \textit{(c,d)} A representative sample under clamped, uniaxial relaxation, respectively, in upper and side views, from which images are extracted for area estimation.}
\end{figure}

Five samples from a single porcine UB were extracted from distinct anatomical locations as shown in Fig.\ref{fig:experimental}(a). The locations are {given and denoted by:} Dorsal (D), lateral (L), lower body (LB), trigone (T), and ventral (V). The samples were {extracted} with a $1 \times 3\, [cm]$ leather punch in the apex-to-base direction as shown in Fig.\ref{fig:experimental} (b). Each sample was clamped, {as illustrated in Figs.\ref{fig:experimental}(c) and (d),} and subjected to five consecutive stress relaxation stages, under prescribed step strains {$\varepsilon_i = 0.25,\,0.50,\,1.00,\,1.50,\,2.00$ over} 30, 45, 45, 45, and 45 minutes, respectively. Besides the strain inputs, the force denoted by {$F^{\text{data}}(t)$} is measured by a $10\,[lb]$ load cell throughout the duration of the test. The cross-sectional area {$A^{\text{data}}(t)$} of each sample is calculated by taking top and side view pictures {that are converted to binary images which are processed in \MATLAB, respectively estimating the base $b(t)$ and height $h(t)$ dimensions of the samples after the application of each consecutive step strain $\varepsilon_i$.} The updated cross-sectional area is assumed to remain constant throughout the relaxation at each strain level, and is evaluated as $A^{\text{data}} = b(t) h(t)$. Given force and area time-series, the true strain is evaluated as $\sigma^{\text{data}}(t) = F^{\text{data}}(t)/A^{\text{data}}(t)$.

\begin{figure}[ht!]
  \includegraphics[width=\columnwidth]{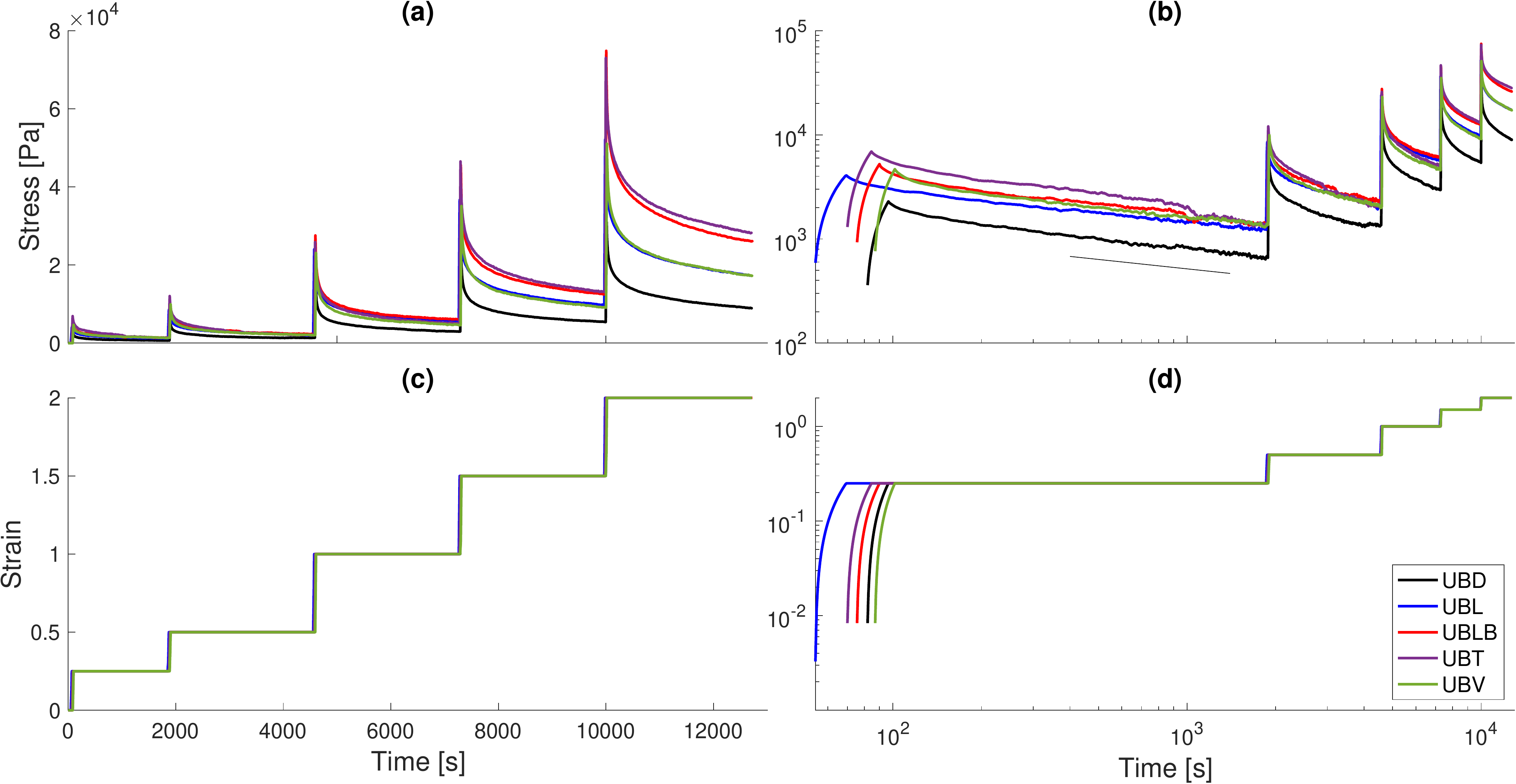}
	\caption{\label{fig:TotalStress}Stress relaxation data for all UB samples. \textit{(a,b)} - stress \textit{vs} time data in linear and logarithmic scales. \textit{(c,d)} - successive step strain \textit{vs} time data in linear and logarithmic scales. The long-term power-law behavior becomes evident especially in the first relaxation step in \textit{(b)}, where the slope of the black line is given by $|\beta| \approx 0.3$.}
\end{figure}

\begin{figure}[ht!]
\centering
  \includegraphics[width=\columnwidth]{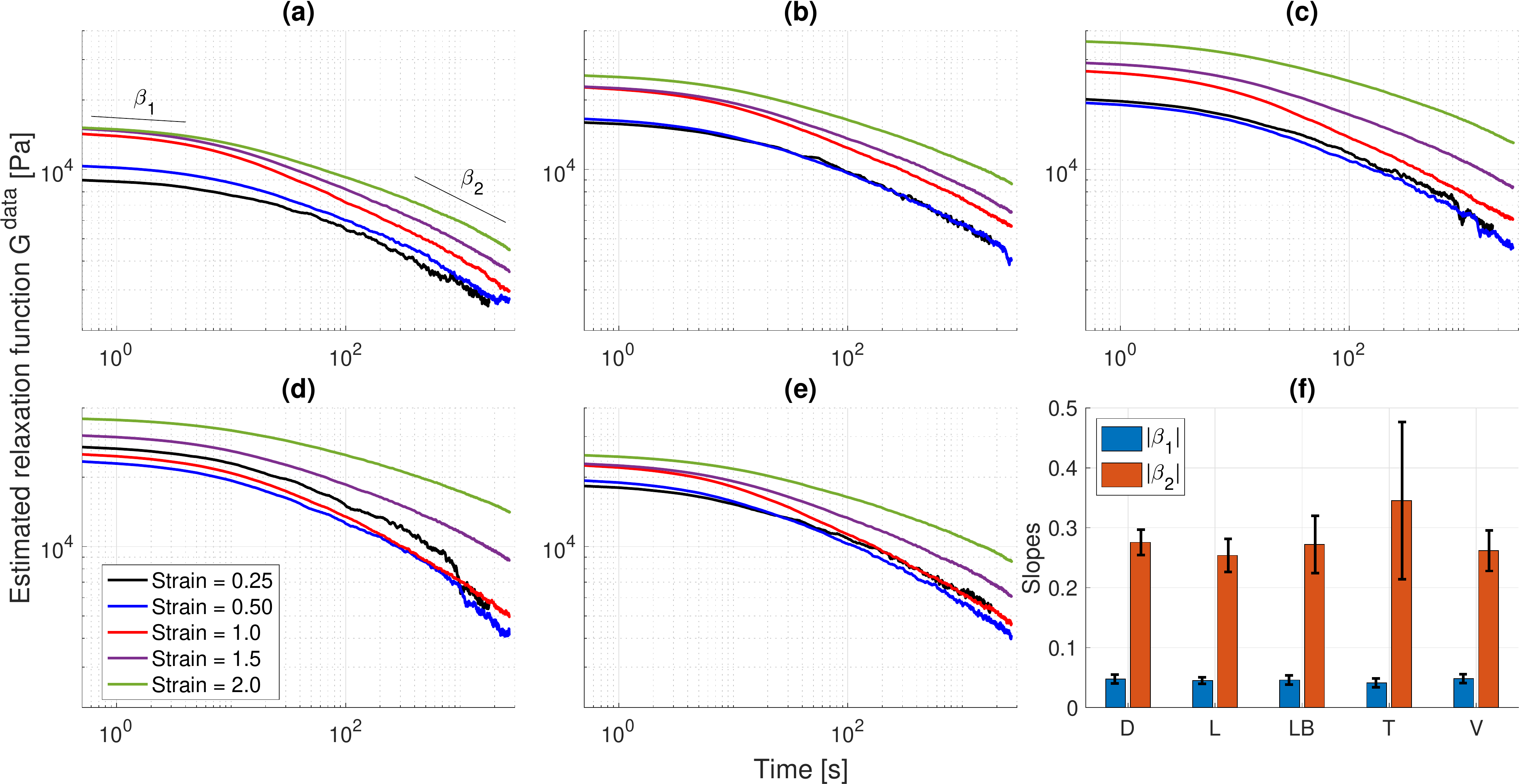}
	\caption{\label{fig:Moduli_all}{Relaxation} moduli for each bladder specimen under distinct applied strains: \textit{(a)} D, \textit{(b)} L, \textit{(c)} LB, \textit{(d)} T, \textit{(e)} V. The relaxation behavior seems to be approximately linear for most specimens in the {$\varepsilon_i \in [0.25, 0.5]$} range. All specimens clearly demonstrate a dependency of both time and applied strains for {$\varepsilon_i > 0.5$}, with largest variations observed in the LB and T samples. Furthermore, the data indicates two distinct relaxation regimes with slower ($t < 3\,[s]$) and faster ($t>400\,[s]$) decays, respectively, with slopes $\beta_1$ and $\beta_2$ illustrated in \textit{(f)}.}
\end{figure}

Once the stress {$\sigma^{\text{data}}(t)$} and strain {$\varepsilon^{\text{data}}(t)$} time-series are obtained, we filter the data through a moving average filter with a time-window of thirty neighbor data points. Figure \ref{fig:TotalStress} illustrates the relaxation curves for all samples in linear and log-log scales. We observe a characteristic power-law scaling for long-time behavior, which is evident in Fig.\ref{fig:TotalStress} (b). As will be shown later through fractional model fits, the relatively low scaling coefficient $\beta$ indicates an anomalous behavior of predominantly elastic nature, {with a \textit{plateau} with low decay rates $\sigma \sim t^{-\beta}$ at larger time-scales (\textit{i.e.}, $t > 400\,[s]$)}. We also note that the trigone and lower body specimens yielded higher stress levels, particularly at very high strains, while the dorsal specimen yielded lower overall values. This is in accordance to stress-strain results obtained by Korossis et al. \cite{korossis2009regional}, that indicated statistically significant, higher collagen phase slopes for the lower body and trigone regions.
We analyze the presence of strain dependency on the UB relaxation behavior, in order to determine if the viscoelasticity is of {linear or nonlinear} nature. Therefore, we employ the definition of linear relaxation modulus $G(t)\,[Pa]$, applied for each fixed strain application from experimental data \cite{haddad1995viscoelasticity}:
\begin{equation}\label{eq:linearrelaxation}
G^{\text{data}}(t) := \sigma^{\text{data}}(t)/\varepsilon_i,
\end{equation}
{with $i=0,\,1,\,\dots,\,N_{\text{steps}}$. In addition to the strain dependency on the relaxation behavior, the above definition also allows us to identify the presence of additional power-laws or stretched exponential behaviors in the tissue response.} Figure \ref{fig:Moduli_all} illustrates the obtained relaxation moduli for all {UB} samples and relaxation steps, {after employing \eqref{eq:linearrelaxation} into the stress time-series data of Fig.\ref{fig:experimental} and performing a translation to a reference initial time, here taken as the first time-step of the data.} We observe that although the relaxation moduli data for each sample {is almost linear for $\varepsilon_i = 0.25,\,0.50$, since the curves approximately overlap (except for the trigone sample)}, the behavior of {$G^{\text{data}}$} is clearly {time- and strain-dependent}. Furthermore, the degree of nonlinearity is more pronounced for the lower body and trigone samples, and less pronounced for the dorsal sample. Interestingly, we notice two limiting power-law behaviors for short and long times. The short time behavior {appears to have a stretched exponential nature, with a limiting power law of smaller magnitude $\beta_1$}, while the long time behavior is associated with a power-law of larger magnitude $\beta_2$, {see Fig.\ref{fig:Moduli_all}(f)}. {Nevertheless, the relaxation behavior clearly transitions from a slower regime to a faster regime, even though the latter still presents a far-from-equilibrum response (no equilibrium glass state) within the experimental time scales.} Larger standard deviations for the long time power-law were observed for the trigone region due to its distinct response for $\varepsilon_0 = 0.25$. We remark that this analysis is just performed to infer the nonlinear relaxation quality the data, and we do not intend to construct a master curve for each {UB} sample.

\subsection{A Fractional Viscoelastic Modeling Framework for Anomalous Tissue Rheology}

{We develop a two-stage framework to assess the qualitative linear/nonlinear relaxation behavior of anomalous materials and identify the most feasible fractional viscoelastic models, in order to alleviate the inherent ill-posedness of model selection problems. In the first stage, we develop an \textit{existence} study to identify the admissible set of anomalous constitutive laws that satisfy the quality of the experimental linear relaxation behavior}, while shedding light on the corresponding microstructural constituents associated to anomalous behavior. {In the second stage, we incorporate large strain effects to the relaxation quality, based on the nonlinear nature of the tissue of interest.}

\subsubsection{{First stage: An existence study of fractional linear viscoelastic models}}
\label{Sec:FLVE}

{Starting with the Scott-Blair (SB) model as the fundamental building block, we construct building block models through parallel and serial combinations to obtain the fractional Kelvin-Voigt (FKV) and fractional Maxwell (FM) models. In our approach, we take into account the anomalous qualities present in the experimental relaxation data and compare them with each of the candidate building block models. Each of the models exhibit distinguished material complexities, such as distinct asymptotic behaviors of relaxation $G(t)$, multiple power-law regimes, slower/faster relaxation at the asymptotic stages \cite{Schiessel1995,Mainardi2010} and presence of material nonlinearities. In the last part of the existence study we classify the candidate models according to their anomalous response, which together with obtained fitting errors and number of material parameters, constitutie the criteria for the model selection procedure. Given the experimental data presented in Section \ref{sec:UBtests}, we focus on the first relaxation step $(\varepsilon_0 = 0.25)$ of Figs.\ref{fig:TotalStress} and \ref{fig:Moduli_all} for the first stage of our data-driven framework}. Our objective is to demonstrate how fractional viscoelastic models are able to capture the UB relaxation with simplistic mechanical arrangements and a small number of material parameters.

\begin{figure}[t]
  \includegraphics[width=\columnwidth]{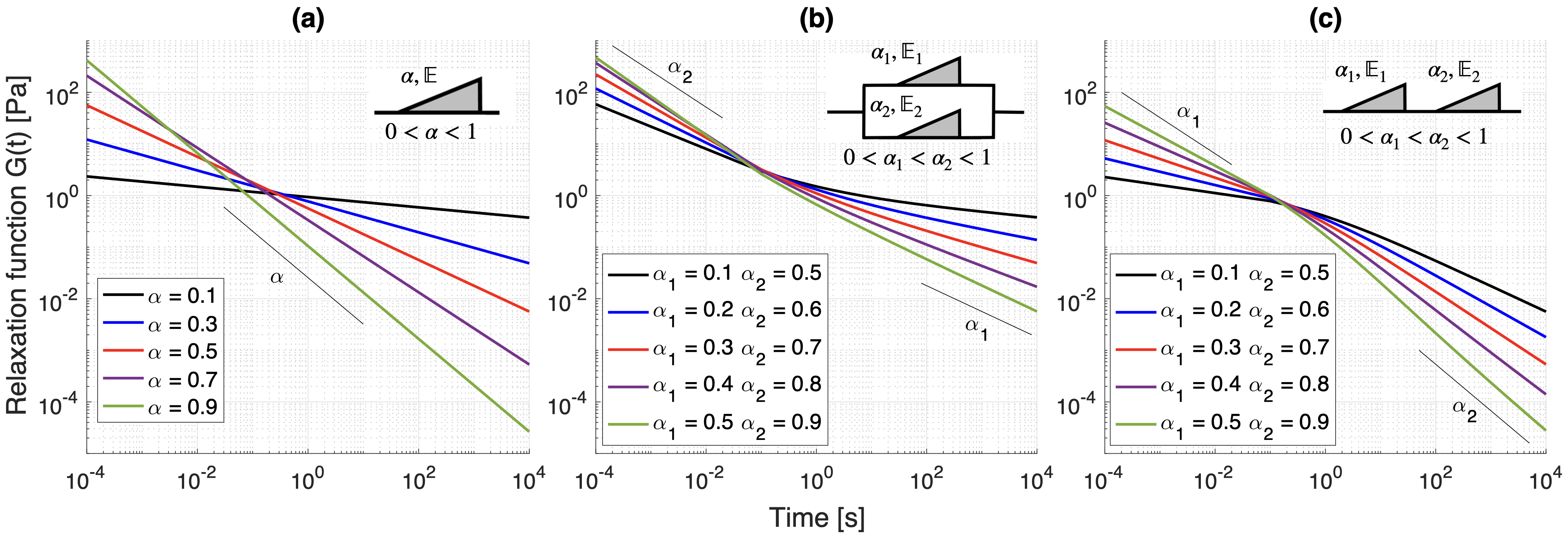}
	\caption{\label{fig:relaxation_models} Relaxation functions $G(t)$ for the building block models under varying fractional models and $\mathbb{E}_1 = 1$, $\mathbb{E}_2 = 1$. \textit{(a)} Scott-Blair, \textit{(b)} Fractional Kelvin-Voigt, and \textit{(c)} Fractional Maxwell. We note the progression from a single, scale-free power-law behavior for the SB model to two dominating power-laws under small and large times for the FKV and FM models.}
\end{figure}

The rheological \textit{building block} for our framework is the fractional {SB} viscoelastic element, which compactly represents an anomalous viscoelastic constitutive law connecting the stresses and strains:
\begin{equation}\label{eq:SB}
\sigma(t) = \mathbb{E} \prescript{\text{C}}{0}{} \mathcal{D} _{t}^{\alpha} \varepsilon(t),
\end{equation}
{with $t > 0$, $\varepsilon(0) = 0$, and} constant fractional order in the range $0 < \alpha < 1$, which provides a material interpolation between Hookean ($\alpha \to 0$) and Newtonian ($\alpha \to 1$) elements. The operator $\prescript{\text{C}}{0}{} \mathcal{D} _{t}^{\alpha} (\cdot)$ represents the time-fractional Caputo derivative given by:
\begin{equation}\label{eq:SBM}
\prescript{\text{C}}{0}{} \mathcal{D} _{t}^{\alpha} u(t) := \frac{1}{\Gamma(1-\alpha)}\int^t_0 \frac{\dot{u}(s)}{(t-s)^\alpha}\,ds,
\end{equation}
where $\Gamma(\cdot)$ represents the Euler-gamma function \cite{Mainardi2010} {and $\dot{u} = du/dt$}. The pair $(\alpha, \mathbb{E})$ uniquely {represents} the SB constants, where the \textit{pseudo-constant} $\mathbb{E}\,[Pa.s^{\alpha}]$ compactly describes textural properties, such as the firmness of the material \cite{Blair1947, Jaishankar2013}. In this sense $\mathbb{E}$ is interpreted as describing a snapshot of a non-equilibrium dynamic process instead of an equilibrium state. The corresponding rheological symbol for the SB model represents a fractal-like arrangement of springs and dashpots \cite{Schiessel1993,McKinley2013}, which we interpret as a compact, upscaled representation of a fractal-like microstructure. Regarding the thermodynamic admissibility {of the SB element and more complex models (i.e., plasticity and damage) involving it}, we refer the reader to Lion \cite{Lion1997} for the SB model, and Suzuki et al. \cite{suzuki2021thermodynamically}. The relaxation function $G^{\text{SB}}(t)\,[Pa]$ for the SB model is given by the following {single,} inverse power-law form:
\begin{equation}\label{eq:SBrelaxation}
G^{\text{SB}}(t) := \frac{\mathbb{E}}{\Gamma(1-\alpha)} t^{-\alpha},
\end{equation}
which is the convolution kernel of the integro-differential form in (\ref{eq:SBM}), {with a modulating pseudo-constant $\mathbb{E}$ for fixed $\alpha$.} Figure \ref{fig:relaxation_models}\textit{(a)} illustrates the behavior of $G^{\text{SB}}(t)$, which is \textit{scale-free}, \textit{i.e.}, a single power-law is present for all $t > 0$. We note that although this relaxation response may seem to be oversimplified, it provides a flexible constitutive interpolation able to, at the very least, take into account the long-term anomalous dynamics of materials, such as the power-law $\beta_2$ in Fig.\ref{fig:Moduli_all}. This also allows the SB element to capture, in {certain} time-scales, power-law behaviors induced by predominantly elastic microstructures, such as collagen networks \cite{doehring2005fractional} with small $\alpha$-values. {We also note that single collagen fibril relaxation data \cite{shen2011viscoelastic} has been calibrated to single-order fractional Kelvin-Zener models (where a single SB element is combined with Hookean springs), with reported values in the range $\alpha = 0.12-0.21$ \cite{bonfanti2020fractional}.} {Similarly, the macroscopic response of ultrasoft tissues such as the brain, which has a more fluid-like character with significant microscopic inter-cell body displacements in the deformation process \cite{reiter2021insights}, has also been captured using SB elements combined with Hookean springs \cite{kohandel2005frequency,davis2006constitutive}. Kohandel et al. \cite{kohandel2005frequency} calibrated a single-order dynamic modulus of a fractional Kelvin-Zener model to bovine brain tissue, obtaining a fractional order $\alpha = 0.6$. Through relaxation experiments, Davis et al. \cite{davis2006constitutive} calibrated a similar model to bovine brain paerenchima also reporting a high value of $\alpha = 0.64$.}

We utilize the SB model as our rheological building block, and define a set of ``building block models", which introduce a higher degree of material complexity through multiple {power-law} behaviors for relaxation and therefore distinct anomalous regimes for small and large time-scales. This multi-fractal type of {rheology} is characteristic of cells \cite{Stamenovic2007} and biological tissues \cite{vincent2012structural}, due to their complex, hierarchical and heterogeneous microstructure. Here we consider {the two simplest} canonical combinations of SB elements. Through a parallel combination, we obtain the fractional Kelvin-Voigt (FKV) model, which has the following stress-strain relationship \cite{Schiessel1993}:
\begin{equation}\label{eq:FKV}
\sigma(t) = \mathbb{E}_1\, \prescript{\text{C}}{0}{} \mathcal{D}^{\alpha_1}_t \varepsilon(t) + \mathbb{E}_2\, \prescript{\text{C}}{0}{} \mathcal{D}^{\alpha_2}_t \varepsilon(t),
\end{equation}
{with $t > 0$, $\varepsilon(0) = 0$, }fractional orders $0 < \alpha_1, \alpha_2 < 1$ and associated pseudo-constants $\mathbb{E}_1\,[Pa.s^{\alpha_1}]$ and $\mathbb{E}_2\,[Pa.s^{\alpha_2}]$. The corresponding relaxation modulus $G(t)$, is also an additive form of {relaxation involving} two SB elements:
\begin{equation}
G^{\text{FKV}}(t) := \frac{\mathbb{E}_1}{\Gamma(1-\alpha_1)} t^{-\alpha_1} + \frac{\mathbb{E}_2}{\Gamma(1-\alpha_2)} t^{-\alpha_2}.
\end{equation}
Figure \ref{fig:relaxation_models}\textit{(b)} illustrates the relaxation function $G^{\text{FKV}}$ for {for several fractional order values}. We notice a response characterized by two power-law regimes, with a transition from faster to slower relaxation slopes. The asymptotic responses for small and large time-scales are given by $G^{\text{FKV}} \sim t^{-\alpha_2}$ as $t\to 0$ and $G^{\text{FKV}} \sim t^{-\alpha_1}$ as $t\to \infty$. We note that this quality allows the FKV model to describe materials that reach an equilibrium behavior for large times when $\alpha_1 \to 0$, which is intuitive from the mechanistic standpoint as one of the SB elements becomes a Hookean spring. {Regarding applications of the FKV model to bio-tissue constituents, Bonfanti et al. \cite{bonfanti2020fractional} has found the combination of low-high fractional orders to be proper for bovine trachea smooth muscle cells, recovering $\alpha_1 = 0.1$ and $\alpha_2 = 0.78$.}

Finally, through a serial combination of SB elements, we obtain the fractional Maxwell (FM) model \cite{Jaishankar2013}, given by:
\begin{equation}\label{eq:FM}
\sigma(t) + \frac{\mathbb{E}_2}{\mathbb{E}_1}\prescript{\text{C}}{0}{}\mathcal{D}^{\alpha_2 - \alpha_1}_t \sigma(t) = \mathbb{E}_2 \prescript{\text{C}}{0}{}\mathcal{D}^{\alpha_2}_t \varepsilon(t),
\end{equation}
{with $t>0$, fractional orders $0 < \alpha_1 < \alpha_2 < 1$, the additional constraint $0 < \alpha_2 - \alpha_1 < 1$,} and two sets of initial conditions for strains $\varepsilon(0) = 0$, and stresses $\sigma(0) = 0$. We note that in the case of non-homogeneous ICs, {one requires compatibility conditions} \cite{Mainardi2010} between stresses and strains at $t=0$. The corresponding relaxation function for this building block model assumes a more complex, Miller-Ross form \cite{Jaishankar2013}:
\begin{equation}\label{eq:G_FMM}
G^{\text{FM}}(t) := \mathbb{E}_1 t^{-\alpha_1} E_{\alpha_2-\alpha_1, 1-\alpha_1}\left(-\frac{\mathbb{E}_1}{\mathbb{E}_2} t^{\alpha_2 - \alpha_1}\right),
\end{equation}
where $E_{a,b}(z)$ denotes the two-parameter Mittag-Leffler function, defined as \cite{Mainardi2010}:
\begin{equation}
  E_{a,b}(z) = \sum^\infty_{k=0} \frac{z^k}{\Gamma(a k + b)},\quad Re(a) > 0,\quad b \in \mathbb{C},\quad z \in \mathbb{C}.
\end{equation}
Interestingly, the presence of a Mittag-Leffler function in (\ref{eq:G_FMM}) {leads} to a stretched exponential relaxation for smaller times and a power-law behavior for longer times, as illustrated in Fig.\ref{fig:relaxation_models}. We also observe that the limit cases are given by $G^{\text{FM}} \sim t^{-\alpha_1}$ as $t\to 0$ and $G^{\text{FM}} \sim t^{-\alpha_2}$ as $t \to \infty$, indicating that the FM model provides a behavior transitioning from slower-to-faster relaxation. {Furthermore, when $\alpha_2 \to 1$, the FM model presents a fluid-like behavior for long times \cite{Schiessel1995}, and therefore allowing a fast relaxation, similar to its integer-order counterpart.} We refer the reader to the {recent work by Bonfanti et al. \cite{bonfanti2020fractional}} for a number of applications of the aforementioned models. We notice that both FKV and FM models are able to recover the SB element with a convenient set of pseudo-constants. {Furthermore, from Fig.\ref{fig:relaxation_models}, we observe that, under a variable order $\alpha(t)$ setting, the SB relaxation function \eqref{eq:SBrelaxation} is able to recover both constant-order FKV and FM models, respectively, through parametric decreasing/increasing variable order functions. However, since variable-order models may potentially increase the material parameter space and computational cost, we restrict our discussion to constant order models.} Furthermore, we also outline more complex building block models that yield more flexible responses, including three to four fractional orders, such as the fractional Kelvin-Zener (FKZ), fractional Poynting-Thomson (FPT), and fractional Burgers (FB), which in turn are able to recover the FKV and FM models. We refer the reader to the works \cite{Schiessel1993,bonfanti2020fractional} for more details on such models.

\subsubsection{{Second stage: Fractional quasi-linear viscoelastic modeling}}

The presented models in Section \ref{Sec:FLVE} provide candidates for power-law relaxation functions that describe the anomalous {linear viscoelasticity} of biotissues, however, {in such materials} the stress-strain relationship {may} becomes nonlinear as {fibers in collagen networks} transition from entangled to aligned with the applied load direction. Therefore, the viscoelastic behavior itself becomes nonlinear and the relaxation function has an intrinsic dependency on the strain levels, as observed in Fig.\ref{fig:Moduli_all} under successive large step-strain applications. To incorporate this additional effect to our modeling framework, we follow \cite{fung2013biomechanics,Craiem2008}, and employ the following quasi-linear, fractional viscoelastic model (FQLV):
\begin{equation}\label{eq:QLV}
		\sigma(t,\varepsilon) = \int^t_0 G(t-s)\frac{\partial \sigma^e(\varepsilon)}{\partial \varepsilon} \dot{\varepsilon}\,ds,
\end{equation}
where the convolution kernel is given by a multiplicative decomposition of a reduced relaxation function $G(t)$ {which may be selected from the first stage of the framework}, and an instantaneous, nonlinear elastic tangent response with stress $\sigma^e$. In the work by Craiem \textit{et al.} \cite{Craiem2008}, the reduced relaxation function has a fractional Kelvin-Voigt-like form with one of the SB {elements} replaced with a Hookean element. Here, we assume a simpler rheology and adopt a Scott-Blair-like reduced relaxation in the form:
\begin{equation}\label{eq:SBGreduced}
  G(t) = E t^{-\alpha} / \Gamma(1-\alpha),
\end{equation}
with the pseudo-constant $E$ with units $[s^\alpha]$, {since the elastic strains $\sigma^e(\varepsilon)$ have units of $[Pa]$}. We adopt the same, two-parameter, exponential nonlinear elastic part as in \cite{Craiem2008}:
\begin{equation}\label{eq:expelastic}
  \sigma^{e}(\varepsilon) = A \left(e^{B \varepsilon} - 1\right),
\end{equation}
with $A$ having units of $[Pa]$ and $B$ being a non-dimensional tuning parameter for the degree of nonlinearity induced by applied strains. Substituting Equations \ref{eq:SBGreduced} and \ref{eq:expelastic} into Eq.\ref{eq:QLV}, we obtain:
\begin{equation}\label{eq:FQLV}
		\sigma(t,\varepsilon) = \frac{E A B}{\Gamma(1-\alpha)}\int^t_0 \frac{e^{B \varepsilon(s)} \dot{\varepsilon}(s)}{(t-s)^\alpha}\,ds,
\end{equation}
which differs slightly from the linear SB model \eqref{eq:SB} in the sense that an additional exponential factor multiplies the function being convoluted. {We remark that reduced relaxation forms for the FKV and FM can also be employed in the FQLV framework. For the FKV model, it would lead to a multi-term convolution of the same nature as \eqref{eq:FQLV} with fractional orders $\alpha_1$ and $\alpha_2$. Therefore, the same numerical methods employed in Section \ref{Sec:numericalmethods} could be employed. However, when dealing with a FM reduced relaxation function, a specialized type of quadrature for a Miller-Ross-type of kernel as \eqref{eq:G_FMM} would be required, potentially with impractical computations of Mittag-Leffler functions for a large number of time-steps.}

\subsubsection{Numerical discretizations}
\label{Sec:numericalmethods}

We discretize the fractional Caputo derivatives in Equations \ref{eq:SB}-\ref{eq:FM} through an implicit L1 finite-difference scheme \cite{lin2007finite}. Therefore, we consider a uniform time-grid with $N$ time-steps of size $\Delta t$, such that $t_n = n \Delta t$, with $n = 0,\,1,\,\dots,\,N$. We remark that although the equations for each of the building block models could be discretized utilizing fast schemes and approaches that treat initial time singularities \textit{(see \cite{zeng2018stable,zhou2020implicit} and references therein)}, the number of utilized data-points is not large, with $N^{data} \approx 3\,000$ for the first relaxation and $N^{data} \approx 25\,000$ for all steps. Also, we avoid the singularity nearby $t \approx 0$ since the first relaxation step is applied at approximately 80 seconds. Nevertheless, the non-smooth nature of the loading would degenerate most of the existing numerical methods for FDEs, and we found that the employed method in this work with the $\Delta t$ described in Section \ref{sec:results} is sufficient for the accuracy to reach the plateau of the experimental data, such that model error is dominant.

In the following, we present the discretized forms for each of our employed linear fractional viscoelastic models. We start with the stresses for the SB model evaluated at $t = t_{n+1}$:
\begin{equation}\label{eq:SBdiscrete}
  \sigma_{n+1} = C_1 \left[ \varepsilon_{n+1} - \varepsilon_{n} + \mathcal{H}^{\alpha_1} \varepsilon \right],
\end{equation}
with {discretization} constant $C_1 = \mathbb{E}/(\Gamma(2-\alpha)\Delta t^{\alpha})$.

{For the FKV model, we obtain:}
\begin{equation}\label{eq:FKVdiscrete}
  \sigma_{n+1} = C_1 \left[ \varepsilon_{n+1} - \varepsilon_{n} + \mathcal{H}^{\alpha_1} \varepsilon \right] + C_2 \left[ \varepsilon_{n+1} - \varepsilon_{n} + \mathcal{H}^{\alpha_2} \varepsilon \right],
\end{equation}
with {discretization} constants $C_1 = \mathbb{E}_1/(\Gamma(2-\alpha_1)\Delta t^{\alpha_1})$ and {$C_2 = \mathbb{E}_2/(\Gamma(2-\alpha_2)\Delta t^{\alpha_2})$}.

{Finally, for the FM model, we have:}
\begin{equation}\label{eq:FMdiscrete}
  \sigma_{n+1} = \frac{C_1 \left[ \varepsilon_{n+1} - \varepsilon_{n} + \mathcal{H}^{\alpha_1} \varepsilon \right] + C_2 \left[ \sigma_{n} - \mathcal{H}^{\alpha_1-\alpha_2} \sigma \right]}{ 1 + C_2 },
\end{equation}
with {discretization} constants $C_1 = \mathbb{E}_1/(\Gamma(2-\alpha_1)\Delta t^{\alpha_1})$ and $C_2 = (\mathbb{E}_1 / \mathbb{E}_2)/(\Gamma(2-\alpha_1+\alpha_2)\Delta t^{\alpha_1-\alpha_2})$. The history terms $\mathcal{H}^\nu u$ {with fractional order $\nu$} in the above equations are given by the following form:
\begin{equation*}
\mathcal{H}^\nu u = \sum^n_{j=1} b_j \left[ u_{n+1-j} - u_{n-j}\right],
\end{equation*}
with weights {$b^\nu_j := (j+1)^{1-\nu} - j^{1-\nu}$}.

The discretization for the FQLV model \eqref{eq:FQLV} employed in this work is shown in \cite{suzuki2021returnmapping}, which is a straightforward, fully-implicit, L1 finite-difference approach with a trapezoidal rule employed on the additional exponential factor. Therefore, the discretized stresses for the FQLV model are given by:
\begin{equation}\label{eq:FQLVdiscrete}
  \sigma_{n+1} = C_1\left[\exp(B \varepsilon_{n+\frac{1}{2}}) \left(\varepsilon_{n+1} - \varepsilon_n\right) + \mathcal{H}^\alpha \left(\varepsilon,\frac{\partial \sigma^e}{\partial\varepsilon}\right)\right],
\end{equation}
with constant $C_1 = E A B / (\Delta t^\alpha \Gamma(2-\alpha))$. The discretized history load in this case is given by:
\begin{equation}
  \mathcal{H}^\alpha \left(\varepsilon,\frac{\partial \sigma^e}{\partial\varepsilon}\right) = \sum^n_{k=1} \exp(B \varepsilon_{n-k+\frac{1}{2}}) \left(\varepsilon_{n-k+1}-\varepsilon_{n-k}\right) b_k,
\end{equation}
with {discretization weights $b^\alpha_k = (k+1)^{1-\alpha} - k^{1-\alpha}$} and $\varepsilon_{i+\frac{1}{2}} = (\varepsilon_i + \varepsilon_{i+1})/2$. The presented discretization has an accuracy of $\mathcal{O}(\Delta t^{2-\alpha})$, and we refer the reader to \cite{suzuki2021returnmapping} for simulations of numerical convergence.

\subsection{Model optimization}

We perform the model fits through a particle-swarm optimization (PSO) algorithm \cite{kennedy1995particle}, which was implemented in \textsf{MATLAB}. The adopted PSO parameters are a population $N_{pop} = 30$ and $N_{it} = 1000$ iterations for the linear cases and $N_{it} = 100$ iterations for the nonlinear cases.
For the linear viscoelastic fits, we set the initial material pseudo-parameter in the $0 \le \mathbb{E} \le 10^8\,[Pa.s^\alpha]$, and the fractional orders are constrained in the $0.0001 \le \alpha \le 0.9999$ range, to ensure that the employed fractional models are able to recover simpler fractional counterparts and also standard rheological elements, if required by the experimental data. For nonlinear cases, we estimate ranges for parameters $A$ and $B$ of the FQLV model by fitting the instantaneous stress response \eqref{eq:expelastic} to each stress peak in every step-strain application of Fig.\ref{fig:TotalStress}. From this preliminary estimate, we have obtained parameters in the ranges $10^4 \le A \le 10^5\,[Pa]$ and $0 \le B \le 2$, which are taken as input parameter ranges for the PSO algorithm. For the relaxation parameters of the FQLV, as in \cite{Craiem2008}, we note that the nature of the power-law relaxation kernel, it is nontrivial to obtain a normalized $G(0^+) = 1$. Nevertheless, for the pseudo-constant we set the range $0 \le E \le 1$ and for the fractional order $\alpha$ we employ the same range as the linear case.

Since the stresses $\sigma^{data}(i)$ and strains $\varepsilon^{data}(i)$ from the relaxation dataset are non-uniform in time, we perform a linear (first-order accurate) interpolation of the strains $\varepsilon^{data}(i)$ to an uniform grid. We then utilize the input strains and compute the global best solution for stress for every PSO iteration through \eqref{eq:SBdiscrete}, \eqref{eq:FKVdiscrete}, \eqref{eq:FMdiscrete}, or \eqref{eq:FQLVdiscrete}. Then, we linearly interpolate the stress back to the nonuniform grid to obtain $\sigma^{model}$. The time-step size for the uniform grid solution is set to $\Delta t = 0.495\, [s]$, which is the minimum time interval between two consecutive data-points. For verification purposes, we tested smaller step-sizes ($\Delta t = 0.0495$) and did not obtain improved results, and we note that all employed numerical discretizations for fractional models are fully-implicit. Therefore, our verification step indicates that model error dominates over discretization error for the employed time-step size. The cost function is defined as:
\begin{equation}\label{eq:Cost}
Cost := \sum^{N_{data}}_{i=0} \left( \sigma^{data}_i - \sigma^{model}_i \right)^2.
\end{equation}

The adopted error measures in this work are the normalized least-squares error (LSE) and root mean squared error (RMSE) between the experimental and mapped simulated stresses, which are respectively given by:
\begin{align*}
LSE & := \frac{\sqrt{\sum^{N_{data}}_{i=1} \left( \sigma^{data}_i - \sigma^{model}_i \right)^2}}{\sqrt{\sum^{N_{data}}_{i=1} \left( \sigma^{data}_i \right)^2}} \times 100\%, \\
RMSE & := \frac{1}{\max(\sigma^{data})}\sqrt{\frac{\sum^{N_{data}}_{i=1} \left( \sigma^{data}_i - \sigma^{model}_i \right)^2}{N_{data}}} \times 100\%. \\
\end{align*}
Finally, all numerical simulations were run in a computer system with Intel Xeon Gold 6148 CPUs with 2.40GHz.

\section{Results and Discussion} \label{sec:results}

\subsection{Linear Viscoelasticity}

\begin{figure}[!ht]
    \subfloat[Scott-Blair model.]{\label{fig:SBM_linear}\includegraphics[width=0.49\columnwidth]{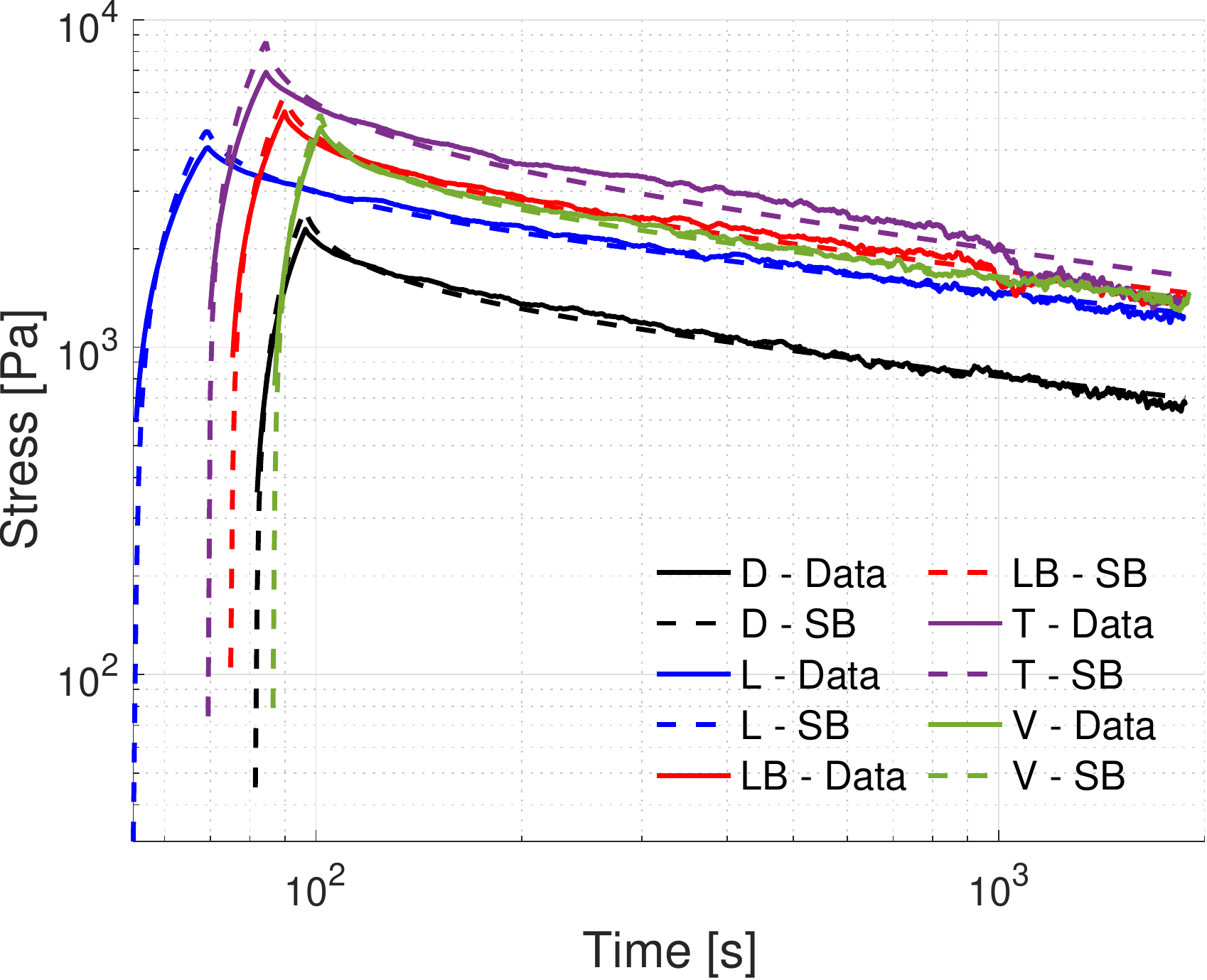}}
    ~
     \subfloat[Fractional Maxwell model.]{\label{fig:FMM_linear}\includegraphics[width=0.49\columnwidth]{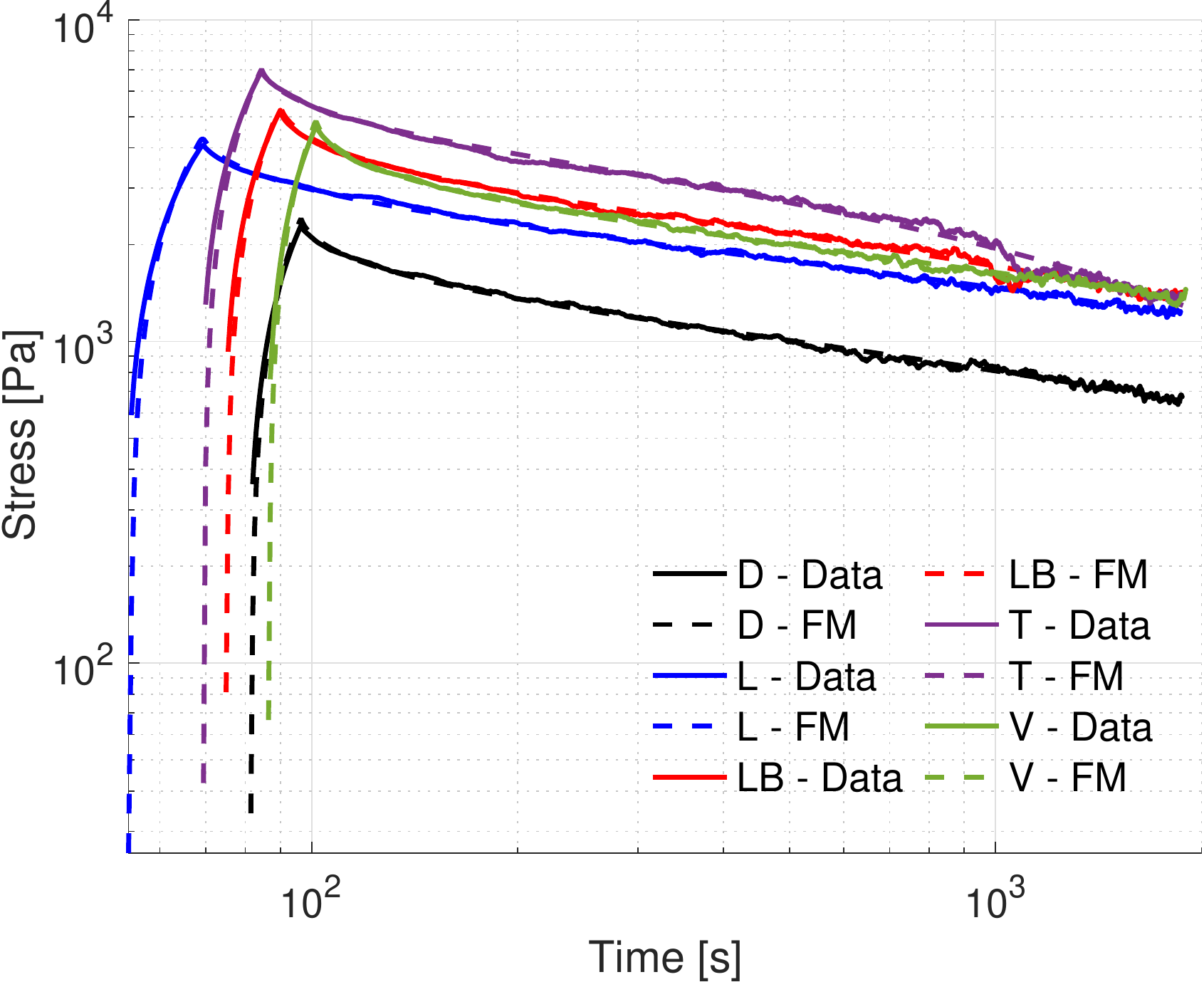}}
     \caption{Obtained linear viscoelastic fits for \textit{(a)} the fractional SB model, and \textit{(b)} fractional Maxwell model for all bladder samples and the first strain step $(\varepsilon_0 = 0.25)$.}
\end{figure}

Figure \ref{fig:SBM_linear} illustrates the obtained fits for all bladder samples utilizing an SB model, for the first strain step $(\varepsilon_0 = 0.25)$. We observe very good fits for most samples, especially at larger time scales, with an exception for the trigone (T) sample due to a sudden stress drop in the experimental data. The fitting quality decreases for all samples at the early relaxation dynamics (nearby the step-strain application), with the SB model underestimating the maximum values of stress peaks. The obtained fractional orders lie in the $0.2-0.3$ range (\textit{see Table \ref{tab:linearfits}}), which are similar to the observed long-term power law from the estimated experimental relaxation functions in Fig.\ref{fig:Moduli_all}. Furthermore, the least-squared errors lie in the $2\% - 11\%$ range, while the RMS errors are within the $2\%-4\%$ range. The higher values of the pseudo-constant and fractional order for the trigone sample are likely due to the SB model accounting for both instantaneous and long-term response over its limited set of two parameters. We also note that the FKV model pragmatically recovered the SB model in all instances, where the PSO algorithm obtained optimal values for the fractional orders that are close to the SB model. In addition, the optimal values for one pseudo-constant is either set to zero, or the summation of $\mathbb{E}_1$ and $\mathbb{E}_2$ recovers the value of $\mathbb{E}$ for the SB model. This indicates the a FKV model does not improve the bladder fit quality, and one would rather employ a SB model with half the amount of material parameters under the same error levels.

Figure \ref{fig:FMM_linear} illustrates the obtained fits for the FM model, where the added flexibility of the underlying power-law/Mittag-Leffler relaxation response improves the fitting quality for both short and long time-scales, yielding least squares errors as low as $2.08\,\%$. The obtained parameters in Table \ref{tab:linearfits} indicate the presence of a predominantly elastic power-law $\alpha_1$ in the $0.17-0.19$ range, and a predominantly viscous $\alpha_2$ in the $0.74-0.99$ range. Particularly, the FM model fit for the trigone specimen indicates the recovery of a dashpot element, and thus the corresponding SB element could be replace by a Newton element. Regarding pseudo-constant values, we note that $\mathbb{E}_1$ values have variations that qualitatively agree with the intensity of stress peaks, but $\mathbb{E}_2$ values can vary in several orders of magnitude, which could be due to the presence of multiple local minima or the discrepancy between obtained fractional orders $\alpha_2$. In general, the dorsal and ventral samples seem to be the most anomalous, as they present both fractional order values sufficiently far from standard elements.

\begin{table}[t]
\centering
\caption{\label{tab:linearfits}Obtained material parameters for all employed linear fractional models and UB samples.}
\begin{tabular}{@{}lcccclccc@{}}
\toprule
\multirow{2}{*}{Model} & \multicolumn{4}{c}{Parameters}    & {} & \multicolumn{2}{c}{Error \%} & \multirow{2}{*}{Sample} \\ \cmidrule(lr){2-5} \cmidrule(lr){7-8}
  & $\mathbb{E}_1[kPa.s^{\alpha_1}]$      & $\alpha_1$ &  $\mathbb{E}_2[kPa.s^{\alpha_2}]$ & $\alpha_2$ & {} &  $LS$ & $RMS$ &    \\ \midrule
SB & 18.1901 & 0.226 & -- & -- & {} & 4.32 & 1.72 & \multirow{3}{*}{D}     \\
FKV& 15.6574 & 0.225 & 2.42019 & 0.229 & {} & 4.32 & 1.72 &   \\
FM& 14.7616 & 0.171 &  3976.90 & 0.743 & {} & 2.29 & 0.91 &    \\
\midrule
SB & 31.3077 & 0.219 &  -- & -- & {} & 3.47 & 1.42 & \multirow{3}{*}{L}     \\
FKV& 31.3462 & 0.220 &  0 & 0.503 & {} & 3.47 & 1.41 &  \\
FM& 26.9311 & 0.186 &  48826.1 & 0.932 & {} & 2.08 & 0.82 & \\
\midrule
SB & 41.6335 & 0.236 &  -- & -- & {}  & 5.33 & 1.98 & \multirow{3}{*}{LB}   \\
FKV& 33.2005 & 0.232 &  7.95429 & 0.252 & {} & 5.34 & 1.99 &    \\
FM& 33.0853 & 0.183 &  37339.9  & 0.935 & {} & 3.23 & 1.20 &    \\
\midrule
SB & 66.7689 & 0.278  &  -- & -- & {} & 11.1 & 3.82 & \multirow{3}{*}{T}     \\
FKV& 66.4714 & 0.278  &  0  & 0.907  & {} & 11.1 & 3.83 &                 \\
FM& 42.4338  & 0.170  &  36938.7 & 0.999  & {} & 4.05 & 1.37 &            \\
\midrule
SB & 34.9254 & 0.220 &  -- & -- & {} & 3.21 & 1.24 & \multirow{3}{*}{V}     \\
FKV& 34.9254 & 0.220 &  0 & 0.579  & {} & 3.21 & 1.25 &                     \\
FM& 30.7605 & 0.188 &  19799.5 & 0.797  & {} & 2.19 & 0.84 &                \\
\bottomrule
\end{tabular}
\end{table}

\begin{figure}[t]
  \includegraphics[width=\columnwidth]{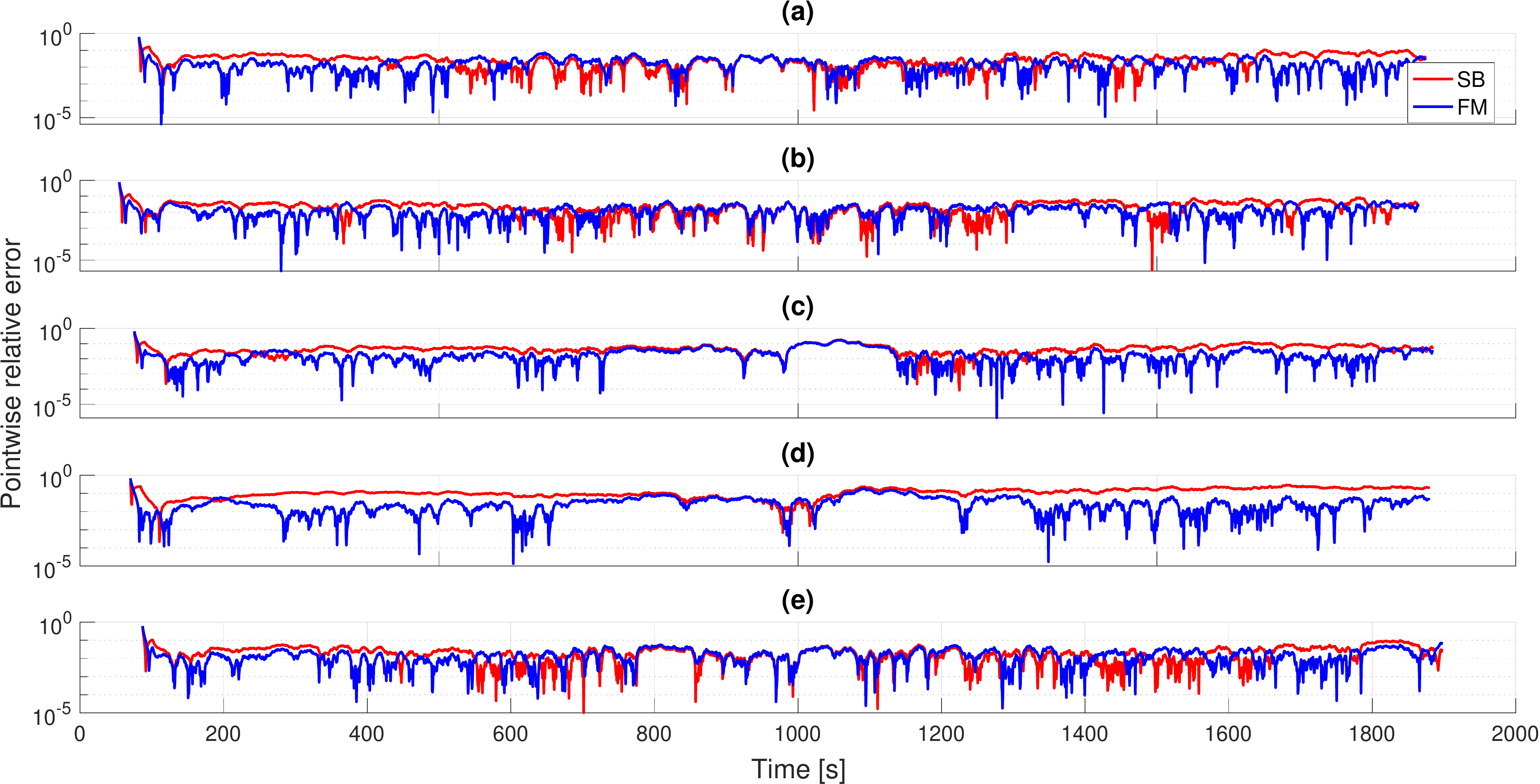}
	\caption{\label{fig:LV_error}Obtained pointwise errors for the linear SB and FM models for the first strain step: \textit{(a)} D, \textit{(b))} L, \textit{(c)} LB, \textit{(d)} T, \textit{(e)} V.}
\end{figure}

Figure \ref{fig:LV_error} illustrates the pointwise relative errors between the SB and FM models and the experimental data for the first strain step. We notice that the SB element has similar errors as the FM model in the $800 < t < 1200\,[s]$ range and larger errors for longer times for most samples. For shorter time ranges, the SB model has larger errors (up to 1 order of magnitude) for all samples. This reinforces the fact that the FM model is more descriptive of both early and long-term dynamics of bladder relaxation, as the qualitative analysis and estimated experimental relaxation moduli suggest. Furthermore, we also note that the better performance of the FM model is also attributed to better approximating the loading ramp and the peak stress preceding the relaxation behavior.

We also employed more complex fractional linear viscoelastic models, such as fractional Kelvin-Zener, Poynting-Thomson and Burgers' (\textit{see \cite{suzuki2021returnmapping} for the models and their corresponding discretizations}). From our employed fitting procedure, all models either recovered or had the same performance as the FM model.

\subsection{Nonlinear Viscoelasticity}

\begin{figure}[!thb]
  \includegraphics[width=\columnwidth]{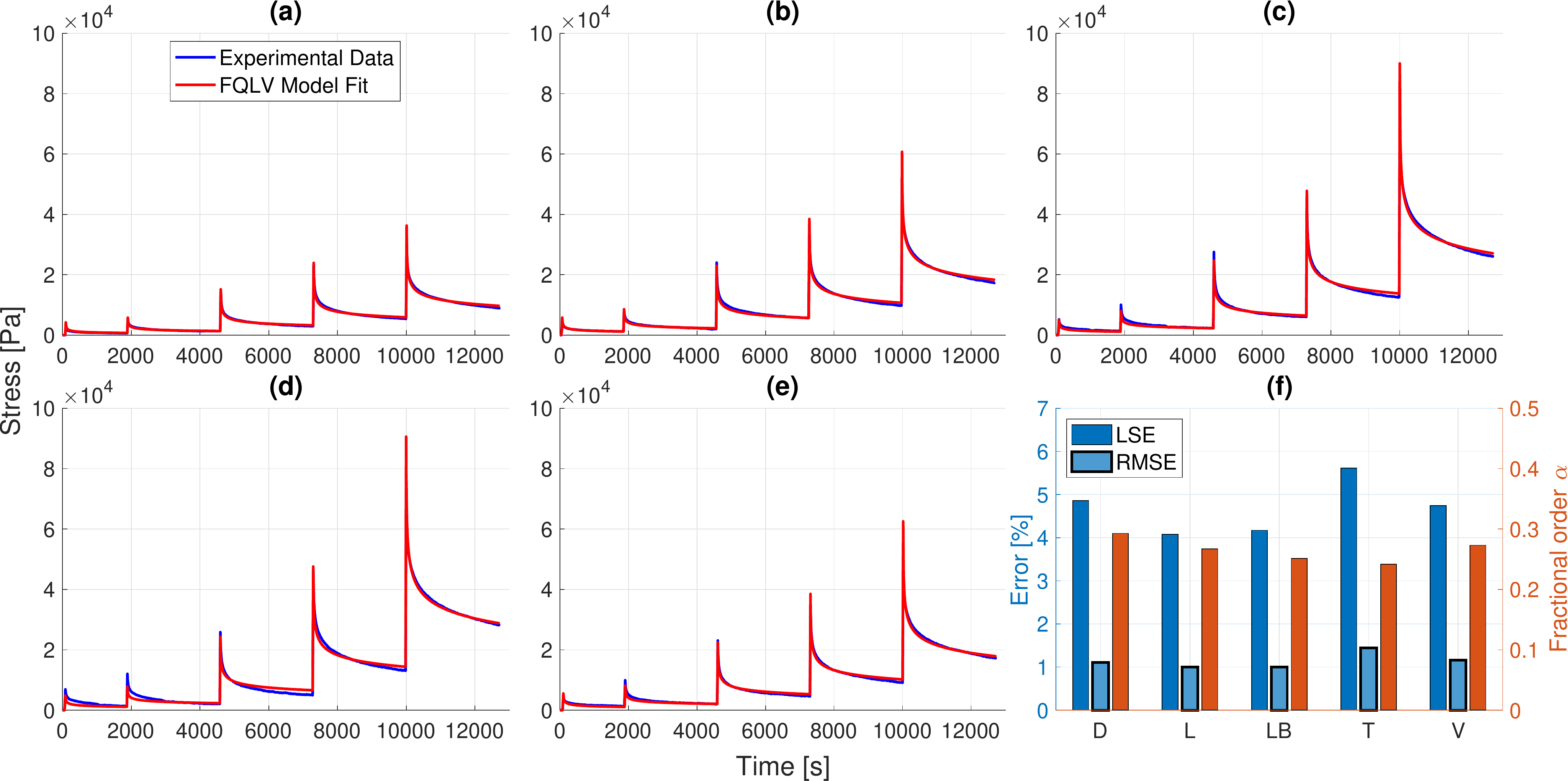}
	\caption{\label{fig:QLV_allsteps}Obtained fits for the fractional QLV model under all strain steps: \textit{(a)} D, \textit{(b)} L, \textit{(c)} LB, \textit{(d)} T, \textit{(e)} V. The fit quality was very good for all bladder samples, with more significant deviations occuring on strain steps 4 and 5. The recovered fractional-orders in \textit{(f)} are within the $0.2 < \alpha < 0.3$ range, which is in accordance with the \textit{a-priori} power-laws obtained from the relaxation moduli data in Fig.\ref{fig:Moduli_all}.}
\end{figure}

Figure \ref{fig:QLV_allsteps} illustrates the obtained fits for the fractional QLV model under all consecutive strain steps, where we observe a very good agreement with the experimental data. Except for the trigone sample, all cases had higher deviations towards the final strain steps. Nevertheless, we note that the error levels are below $6\%$ (LSE) and $2\%$ (RMSE) for the entire dataset, under 4 material parameters, which are listed in Table \ref{tab:FQLVfits}. Furthermore, the obtained fractional-orders lie in the range $0.24-0.3$ which are in accordance with the estimated power-laws in our \textit{a-priori} analysis presented in Fig.\ref{fig:Moduli_all}. Particularly, the lowest fractional order was obtained for the trigone specimen, and highest for the dorsal one. A slightly higher degree of nonlinearity is also recovered for the trigone and lower-body samples due to the larger values of $B$.

\begin{table}[!thb]
\centering
\caption{\label{tab:FQLVfits}Obtained material parameters for the FQLV model with all UB samples.}
\begin{tabular}{@{}llcclcclrr@{}}
\toprule
\multirow{2}{*}{Sample} & {} & \multicolumn{5}{c}{Parameters} & {} & \multicolumn{2}{c}{Error \%} \\ \cmidrule(lr){2-7} \cmidrule(lr){9-10}
  & {}  & $A\,[kPa]$      & $B$ & {} & $E\,[s^{\alpha}]$ & $\alpha$ & {} & $LSE$ & $RMSE$                      \\ \midrule
D & {}  & 53.8823 & 0.7803 & {} & 0.7298 & 0.2928 & {} & 4.85 &  1.11 \\
L & {}  & 79.1646 & 0.8823 & {} & 0.5677 & 0.2673 & {} & 4.08 &  1.00 \\
LB & {} & 74.5369 & 1.2192 & {} & 0.3463 & 0.2510 & {} & 4.17 &  1.00 \\
T & {}  & 63.4435 & 1.2642 & {} & 0.3590 & 0.2419 & {} & 5.61 &  1.44 \\
V & {}  & 59.3282 & 0.9449 & {} & 0.6704 & 0.2732 & {} & 4.74 &  1.16 \\
\bottomrule
\end{tabular}
\end{table}

Finally, Figure \ref{fig:QLV_error_allsteps} illustrates the pointwise relative errors for the FQLV model under all bladder samples and strain steps. We observe a similar error behavior as the SB model in Fig.\ref{fig:LV_error}, particularly for larger applied strains.

\subsection{Discussion}

To our best understanding, this was the first work in the literature addressing fractional viscoelastic modeling to bladder tissues. From the results obtained for our building block models, the overall lower range of fractional orders obtained for all linear/nonlinear models is $0.17-0.3$, indicating a predominantly elastic yet highly anomalous behavior with smaller decay rates at long times, \textit{i.e.}, the presence of far-from-equilibrium dynamics. A similar parametric range was obtained in other anomalous systems such as arterial wall relaxation \cite{Craiem2008}, aortic valve tissue \cite{doehring2005fractional}, 1D and 3D brain artery walls under fluid-structure interactions \cite{perdikaris2014fractional,yu2016fractional}, canine and bovine liver tissue \cite{kiss2004viscoelastic,capilnasiu2020nonlinear}, and lung tissue \cite{Suki1994}. As suggested by Doehring \textit{et al.} \cite{doehring2005fractional}, small $\alpha$-values can be indications of strong fractality in bio-tissue microstructure such as collagen fibers, which are vastly present in the UB, and particularly with a larger network in the trigone region. The larger values of fractional orders $\alpha_2$ in the $0.74-0.99$ range obtained by the fractional Maxwell model is similar to those obtained for brain tissue relaxation \cite{davis2006constitutive} and human ear \cite{naghib2018ear,naghibolhosseini2015estimation}. This indicates a significantly more dissipative behavior, possibly compensating the highly-anomalous behavior provided by the smaller fractional order $\alpha_1$ for short time-scales, and thus better fitting the slower relaxation nearby the load application. The transitional behavior from slower-to-faster relaxation slopes observed from the UB specimens and captured by the FM model were also noticed in bio-tissues composed of weakly cross-linked collagen networks \cite{vincent2012structural}. We note that although the captured fractional orders $\alpha_1,\,\alpha_2$ for the FM model on the UB relaxation do not quantitatively match the slopes $\beta_1,\,\beta_2$ for $G^{data}(t)$ in Fig.\ref{fig:Moduli_all}, these fractional orders refer to the asymptotic behaviors at $t \to 0$ and $t \to \infty$, as illustrated in our existence study in Fig.\ref{fig:relaxation_models}, and it is likely that  relaxation experiments under a larger range of time-scales would yield a better quantitative agreement. For the purpose of our existence study, we consider a qualitative agreement and small error levels to be sufficient to select a valid candidate building block model.

The nonlinear viscoelastic behavior was well approximated by the employed FQLV model, which decomposes the relaxation kernel in a multiplicative fashion into a power-law reduced relaxation function and a tangent elastic stiffness described by an exponential elastic stress form. This allowed the nonlinear part of our existence study to capture the complex rheology of the UB with large applied strains (up to $200\%$) and RMS errors as low as $1\%$. In fact, Jokandan \textit{et al.} \cite{JOKANDAN201892} observed an exponential-like stress-strain response in quasi-static tensile testing of porcine bladder samples. Specifically, under relaxed states, the entangled configuration of collagen fibers yield a linear stress strain relationship, but a nonlinear regime with much higher stress levels is attained once the fibers align with the load direction and store most of the strain energy in the system. Korossis \textit{et al.} \cite{korossis2009regional} attributed the linear region to be predominantly driven by elastin, and the nonlinear phase by collagen.

\begin{figure}[!thb]
  \includegraphics[width=\columnwidth]{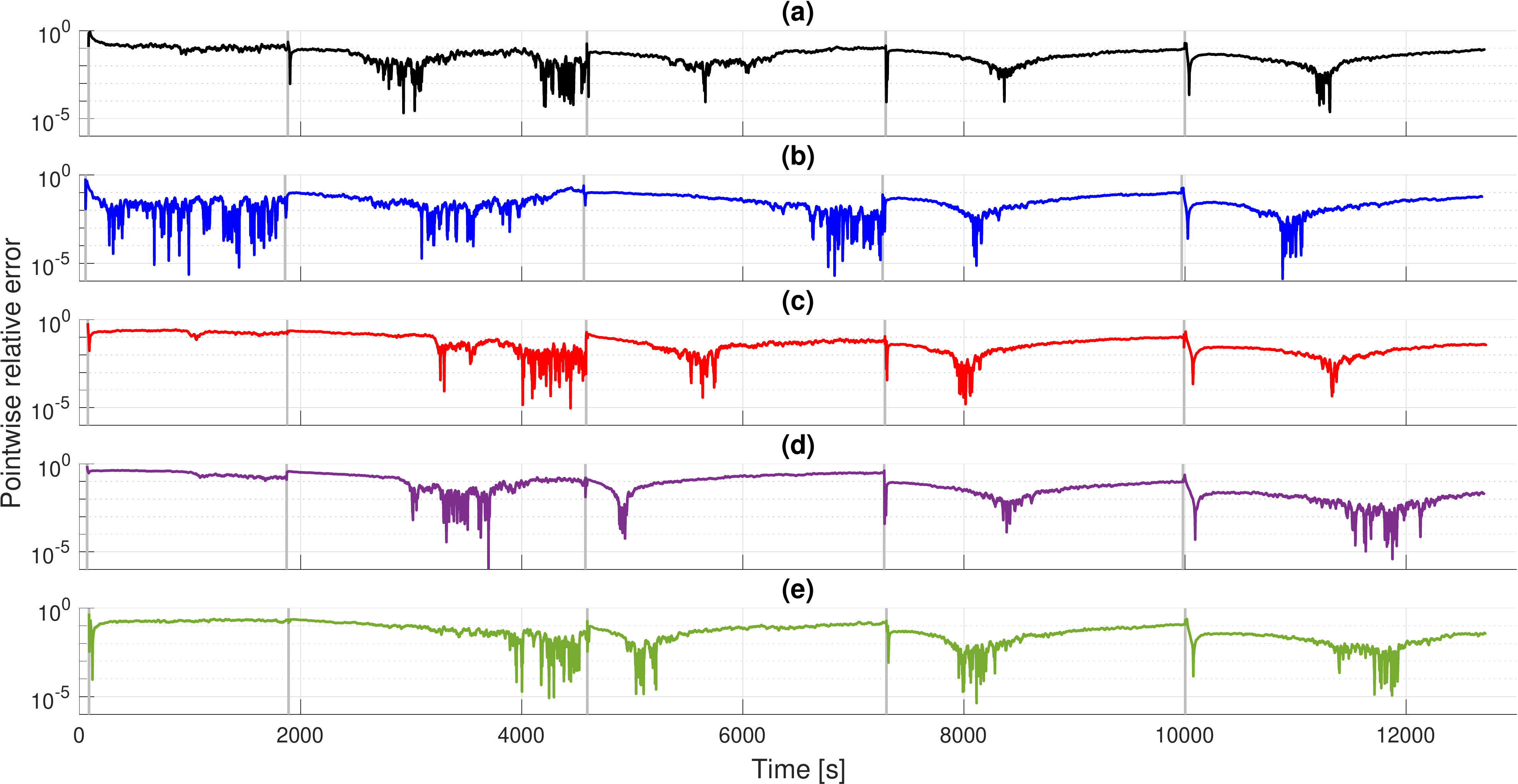}
	\caption{\label{fig:QLV_error_allsteps}Obtained pointwise errors for the FQLV model under all strain steps: \textit{(a)} D, \textit{(b))} L, \textit{(c)} LB, \textit{(d)} T, \textit{(e)} V. The vertical gray lines represent the step strain application instants.}
\end{figure}

Our pointwise relative error analysis reinforced the idea that the FM model is more descriptive of both slower early dynamics of the porcine bladder, but also the faster dynamics observed at longer time-scales in the estimated experimental relaxation modulus. In our error analysis for the nonlinear case, Fig.\ref{fig:QLV_error_allsteps} indicates a similar qualitative error behavior between the FQLV model and the linear SB element in Fig.\ref{fig:LV_error}, especially towards the larger strain regime, with higher errors in the small and large times after the step-strain applications. This suggests that coupling a fractional Maxwell-type reduced relaxation function to the existing framework would very likely improve our results for the nonlinear case as it did in the linear one.

Regarding our developed framework, the existence study proved to be interesting to identify the most proper fractional linear viscoelastic model for stress relaxation, which can later inform the fractional quasi-linear viscoelastic model on the proper form of the reduced relaxation function. For the UB, we conclude that while the SB, FKV and FM models yield errors in the same order of magnitude, the FM model better captures the two power-law qualitative behavior of the data, which is fundamental for both short- and long-term predictions of tissue response. Nevertheless, the SB model provided satisfactory results for the observed experimental time-scale, and the FKV model proved to be redundant and a source of ill-posedness in a model selection framework, since it obtained the same performance as the SB model with twice as many parameters. Although we cannot guarantee that the obtained model parameters provide a global minimum for the cost function (\ref{eq:Cost}), we find our obtained fitting errors, increased number of material parameters, and diverging qualitative behaviors between the experimental data in Fig.\ref{fig:Moduli_all} and the FKV relaxation behavior from Fig.\ref{fig:relaxation_models} to be sufficient to exclude the FKV as a viable candidate for the UB. The same analysis applies to other tested models not shown here, such as fractional Kelvin-Zener, fractional Poynting-Thomson and fractional Burgers' models, which consistently recovered the FM model.

Regarding potential improvements, anisotropy of bladder samples has been observed in existing studies \cite{MORALESORCAJO2018263}, and therefore an interesting step would be to map the variation of fractional-orders for distinct sample orientations. Based on the obtained linear and nonlinear results, another interesting aspect would be to incorporate a fractional Maxwell-type relaxation function to the fractional QLV framework, similar to Doehring \textit{et al.} \cite{doehring2005fractional} to better describe the initial relaxation process in each strain step. However, due to the larger number of time-steps required by our dataset, an efficient numerical method would be required to handle the resulting differ-integral with a Miller-Ross relaxation kernel, or the FQLV framework would need to be developed in differential form, likely as a system of equations comprised of an FDE solving for elastic stresses and a separate equation for nonlinear elasticity. Finally, biaxial tests and models would give insight in the effects of shear stress to the tissue behavior, and confronting creep predictions with experiments would allow one to verify the consistency of the obtained parameters, as already succesfully done with other anomalous materials \cite{Jaishankar2013}.

\section{Conclusions}
\label{Sec:Conclusions}

We developed a data-driven fractional modeling framework for linear and nonlinear viscoelasticity which was validated for the first time in the uniaxial relaxation of porcine urinary bladder tissue for a wide range of applied strains. Our approach employed fractional linear and quasi-linear viscoelastic models to account for anomalous power-law relaxation and large strains. Our main findings in this study were:
\begin{itemize}
\item Our existence study was able to relate the power-law features of specimens from five distinct bladder samples to the fractional orders in linear/quasi-linear fractional models, and is an interesting step towards an automated model selection framework.

\item The bladder uniaxial relaxation data was obtained from consecutive and increasing step strain applications, indicating the presence of nonlinear, strain-dependent effects on the relaxation functions.

\item Our obtained data is consistent with other bladder rheology studies, indicating higher stress levels for the trigone region, and an exponential-like stress-strain relationship.

\item Among linear viscoelastic models employed for the first relaxation step ($25\%$ strains), the fractional Maxwell model was the most suited for all regional bladder samples, with two fractional orders, which dominate short- and long-times. More complex linear fractional models consistently recovered the FM, and the FKV model proved to be unfeasible since it recovers a SB element.

\item The employed fractional quasi-linear viscoelastic model successfully captured the multi-step relaxation behavior with four material parameters and without any requirement of parameter recalibration, yielding a fractional order range $\alpha = 0.25 - 0.30$ with root mean squared errors below $2\%$.

\item Fractional calculus can be a interesting modeling alternative to describe the linear/nonlinear behavior of porcine UB especially within a material model selection framework, since fractional models potentially provide a reduced number of material parameters due to the presence of multiple power-laws in relaxation.
\end{itemize}
Regarding potential future steps, investigating the possibility of plastic deformations under large strains by employing quasi-linear visco-plastic effects \cite{suzuki2021returnmapping,Suzuki2016,suzuki2021thermodynamically} and also failure mechanisms \cite{barros2021integrated} would be interesting studies towards the life-cycle prediction of such anomalous bio-tissues.
Finally the variation of tissue properties by anatomical location, orientation, and layers of the UB motivates further studies on distributed order models incorporating nonlinearities. In this sense, mathematical and computational frameworks employing distributed-order viscoelastic models that learn distributions for fractional orders would be able to account for multi-fractal heterogeneous media and stochastic effects, leading to predictive zero-dimensional models with a reduced number of parameters.

\section*{Acknowledgments}
J.L. Suzuki and M. Zayernouri were supported by the ARO YIP Award (W911NF-19-1-0444), the MURI/ARO (W911NF-15-1-0562) and the NSF Award (DMS-1923201). The computational resources and services were provided by the Institute for Cyber-Enabled Research (ICER) at Michigan State University

\bibliographystyle{siamplain}
\bibliography{references}

\end{document}